\documentclass[12pt]{amsart}
\usepackage[dvipdf]{graphicx}
\usepackage{times}
\usepackage{geometry}
\newtheorem{thm}{Theorem}
\newtheorem{conj}[thm]{Conjecture}
\newtheorem{cor}[thm]{Corollary}
\newtheorem{que}[thm]{Question}

\title[Finite planar emulators for $K_{4,5} - 4K_2$ and
  $K_{1,2,2,2}$]{Finite planar emulators for $K_{4,5} - 4K_2$ and
  $K_{1,2,2,2}$ \\[5pt] and Fellows' Conjecture}

\date{\today}  
\address{Department of mathematical Sciences, University of Arkansas, 
Fayetteville, AR 72701} %
\address{Department of Information and Computer Sciences, Nara Women's
  University  Kitauoya Nishimachi, Nara 630-8506, Japan} %
\email{yoav@uark.edu} %
\email{yamasita@ics.nara-wu.ac.jp} %
\author{Yo'av Rieck} %
\author{Yasushi Yamashita} %


\begin{document}

\begin{abstract}
In 1988 Fellows conjectured that if a finite, connected graph
admits a finite planar emulator, then it admits a finite planar
cover.  We construct a finite planar emulator for $K_{4,5} - 4K_2$.
D. Archdeacon~\cite{archdeacon2} showed that $K_{4,5} - 4K_2$ does not
admit a finite planar cover; thus $K_{4,5} - 4K_2$ provides a
counterexample to Fellows' Conjecture.  

It is known that Negami's Planar Cover Conjecture is true if and
only if $K_{1,2,2,2}$ admits no finite planar cover.  We construct a
finite planar emulator for $K_{1,2,2,2}$.  The existence of a finite
planar cover for $K_{1,2,2,2}$ is still open.
\end{abstract}

\nocite{*}

\maketitle

\section{Introduction}

We begin by defining the main concepts used in this paper. All graphs
considered are assumed to be finite and simple.  A map between graphs
is assumed to map vertices to vertices and edges to edges. Let
$\widetilde{G}$ and $G$ be graphs.  We say that $\widetilde{G}$ is  a
{\it cover} (resp. {\it emulator}) of $G$ if there exists a map
$f:\widetilde{G} \to G$ so  that $f$ is surjective and for any vertex
$\tilde v$ of $\widetilde{G}$, the map induced by $f$ from the
neighbors of $\tilde v$ to the neighbors of $f(\tilde v)$ is a
bijection (resp. surjection).   A cover (resp. emulator) is called
{\it regular} if there is a subgroup $\Gamma \le Aut(\widetilde{G})$ (the
automorphism group of $\widetilde{G}$) so
that $G \cong \widetilde{G}/\Gamma$, and $f$ is equivalent to the
natural projection.  In this paper, regular covers and emulators are
only used when citing results of Negami and Kitakubo; for detailed
definitions see \cite{negami} (for covers) and \cite{kitakubo} (for
emulators).  We note that Kitakubo used the term {\it branched covers}
for emulators.

Let $i:S^2 \to \mathbb{R}P^2$ be the projection from the sphere to the
projective plane given by identifying antipodal points.  If a graph
$G$ embeds in $\mathbb{R}P^2$, then $i^{-1}(G)$ is a
planar double cover of $G$. Conversely, in \cite{negami} Negami
proved that if a connected graph $G$ admits a finite planar {\it
  regular} cover, then $G$ embeds in $\mathbb{R}P^2$.  Negami
conjectured that this holds in general:

\begin{conj}[Negami's Planar Cover Conjecture]
\label{conj:negami}
A connected graph has a finite planar cover if and only if it embeds
in the projective plane. 
\end{conj}

Kitakubo generalized Negami's theorem, showing that if a graph has a
finite planar regular emulator, then it embeds in the projective
plane.  (The authors gave a further generalization in \cite{yr}.)
The following conjecture appears in
\cite[Conjecture~2]{negami-survey}, where Negami attributes it to
Kitakubo: 

\begin{conj}
\label{conj:kitakubo}
A connected graph has a finite planar emulator if and only if it
embeds in the projective plane. 
\end{conj}

Prior to Kitakubo, planar emulators were studied by Fellows
\cite{fellowsphd}\cite{fellows}, who 
posed the conjecture below; see, for example,
\cite[Conjecture~4]{hlineny} or \cite{negami-survey}:

\begin{conj}[Fellows]
\label{conj:fellows}
A connected graph has a finite planar emulator if and only if it has a
finite planar cover.
\end{conj}

In \cite{hlineny} Hlin{\v{e}}n{\'y} constructed a graph that admits an
emulator that embeds in the genus 3 surface, but does not admit a cover
that embeds there.

In this note we prove:

\begin{thm}
\label{thm}
The graphs $K_{4,5} - 4K_2$ and $K_{1,2,2,2}$ admit finite planar
emulators. 
\end{thm}

Archdeacon~\cite{archdeacon2} proved that $K_{4,5}-4K_2$ does not
admit a finite planar cover.  Together with Theorem~\ref{thm}, we
get:  

\begin{cor}
The graph $K_{4,5} - 4K_2$ gives a counterexample to
Conjectures~\ref{conj:kitakubo} and \ref{conj:fellows}.
\end{cor}

It is known that $K_{1,2,2,2}$ does not embed in $\mathbb{R}P^2$ 
\cite{103}.  Hence, if it admits a finite planar cover, 
Negami's Planar Cover Conjecture is false.  The work of
Archdeacon, Fellows, Hlin{\v{e}}n{\'y}, and Negami shows that
the converse also 
holds, and Negami's Planar Cover Conjecture is in fact equivalent to
$K_{1,2,2,2}$ having no finite planar cover; see, for example,
\cite{negami-survey} or \cite{hlineny2} and references therein.   At
the time of writing, 
the existence of a finite planar cover to $K_{1,2,2,2}$ remains an
intriguing open question.  However,
Theorem~\ref{thm} shows that $K_{1,2,2,2}$ does admit a finite planar
emulator.  Perhaps this should not be seen as evidence against
Negami's Planar Cover Conjecture.  Perhaps this should be seen as
evidence that finite planar emulators are ubiquitous (although clearly 
not all graphs have finite planar emulators).  We note that if
Negami's Planar Cover Conjecture holds, then the existence of a finite
planar cover can be decided in linear time~\cite{mohar}; the set of
forbidden minors is given by Archdeacon~\cite{archdeacon-projective}.
By Robertson and Seymour~\cite{rs} there is a set of forbidden minors
for existence of a finite planar emulator; therefore~\cite{rs2}
existence of such an emulator can be decided in polynomial time.

\begin{que}
\label{que}
What graphs admit finite planar emulators?  What is the set of
forbidden minors?  Construct an algorithm to decide if a given graph
admits a finite planar emulator.  What is the complexity of this
problem? 
\end{que}

In Section~\ref{sec:ArchdeaconsGraph} we explicitly show an emulator
with 50 vertices for $K_{4,5} - 4K_2$ and in Section~\ref{sec:1222} we
explicitly show an emulator with 266 vertices for $K_{1,2,2,2}$, thus
proving Theorem~\ref{thm}.  The emulator for $K_{1,2,2,2}$ is
symmetric and quotients out to an emulator with 133 vertices that embeds 
in $\mathbb{R}P^2$.

\medskip

\noindent{\bf Acknowledgment.}  We thank Mike Fellows and Petr
Hlin{\v{e}}n{\'y} for  
helpful comments regarding Question~\ref{que},  
and the referees for their work and remarks.
The first named author
thanks  Tsuyoshi Kobayashi and the department of mathematics of Nara
Women's University for their hospitality during the time this research
was conducted.

\section{A Finite Planar Emulator for $K_{4,5} - 4K_2$.}
\label{sec:ArchdeaconsGraph}

\begin{figure}
\begin{center}
\includegraphics{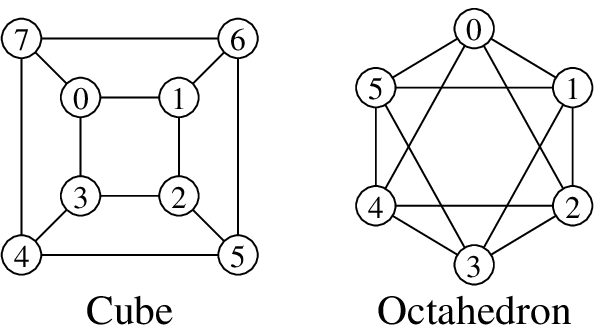}
\end{center}
\caption{}
\label{fig:cube}
\end{figure}

\begin{figure}
\begin{center}
\includegraphics{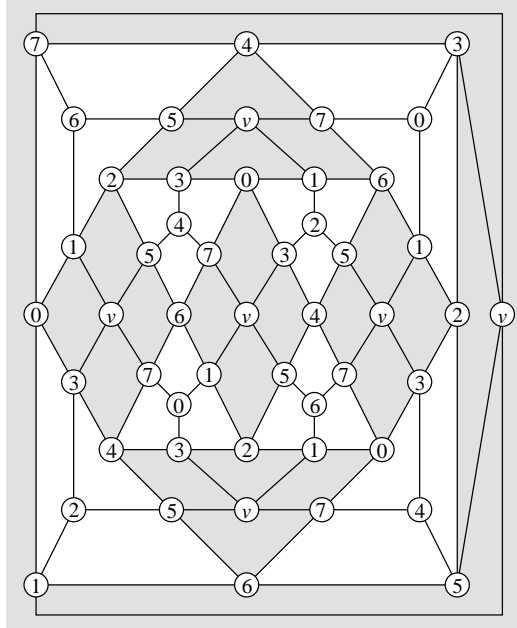}
\end{center}
\caption{A finite planar emulator of $K_{4,5} - 4 K_2$}
\label{fig:ArchdeaconsGraph}
\end{figure}

An emulator for $K_{4,5} - 4K_2$ is given in
Figure~\ref{fig:ArchdeaconsGraph}.  We explain how to read this
graph.  $K_{4,5} - 4K_2$ is constructed as follows: start with the 
1-skeleton of a cube (Figure~\ref{fig:cube}) and add a ninth vertex,
denoted $v$, that is connected to the vertices of the cube labeled 1,
3, 5, and 7.  The graph shown in Figure~\ref{fig:ArchdeaconsGraph}
maps to $K_{4,5} - 4K_2$; each vertex is shown as a small circle
labeled by the vertex of $K_{4,5} - 4K_2$ it gets sent to.  It can be
checked directly that it is a finite planar emulator of $K_{4,5} -
4K_2$.

\medskip

\noindent{\it Remarks.} 
\begin{enumerate}
\item It is easy to see how the graph was put together.  It is made of 8
white ``triangles'', each triangle meeting 3 others (this pattern can
be seen by taking the convex hull of the midpoints of the edges of a
cube or an octahedron).  Each triangle is simply a corner of the cube
(3 squares).
\item Note that the graph presented in Figure~\ref{fig:ArchdeaconsGraph}
is not a cover of $K_{4,5} - 4K_2$.  For example, we can find a
vertex with label 0 which is adjacent to two vertices labeled 3.
\end{enumerate}

\section{A Finite Planar Emulator for $K_{1,2,2,2}$.}
\label{sec:1222}

\begin{figure}
\begin{center}
\includegraphics[width=145mm]{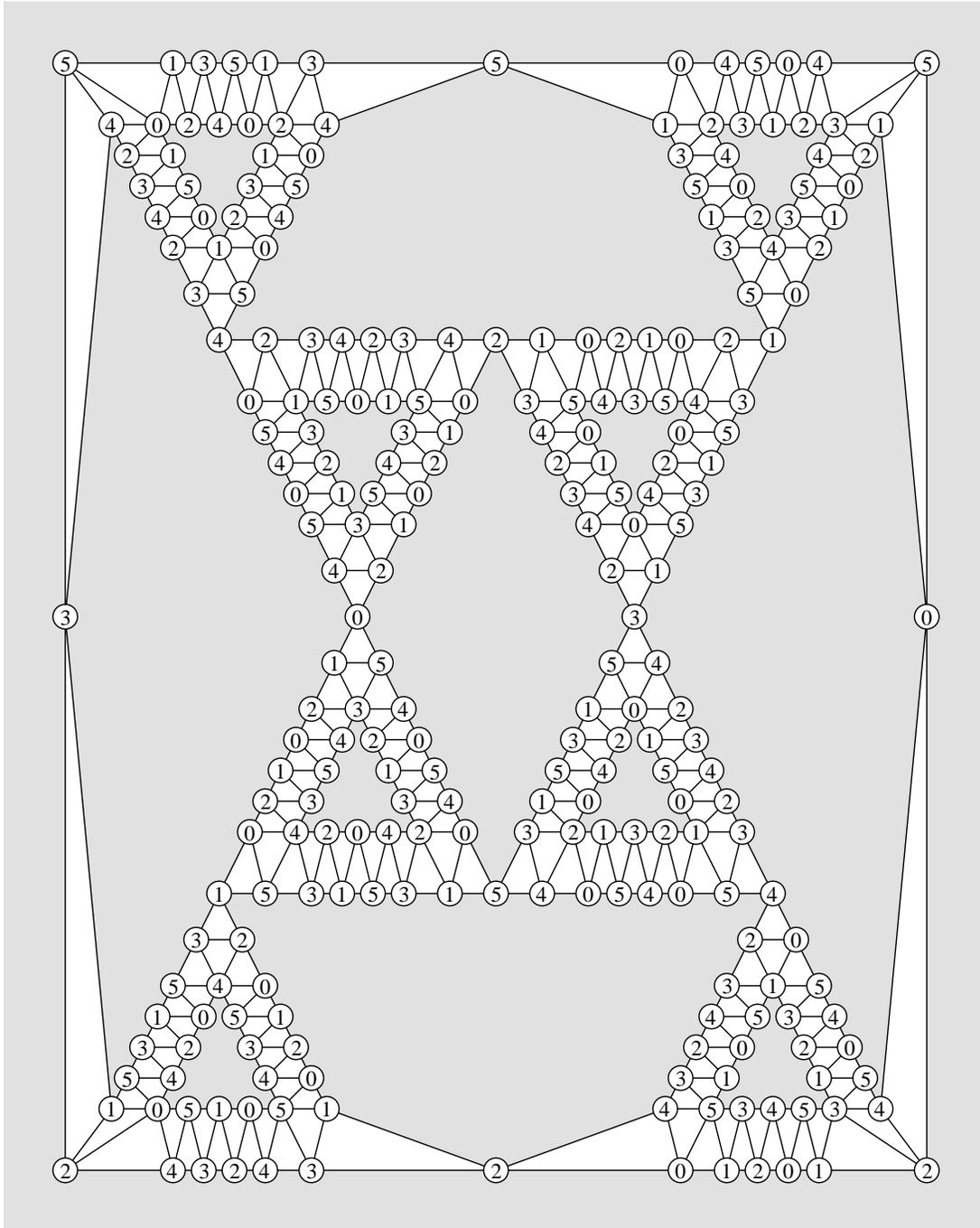}
\end{center}
\caption{A finite planar emulator of $K_{1,2,2,2}$}
\label{fig:1222}
\end{figure}

An emulator for $K_{1,2,2,2}$ is derived from Figure~\ref{fig:1222}.  
We explain how to read this graph.  $K_{1,2,2,2}$ is constructed as
follows: start with the 1-skeleton of an octahedron
(Figure~\ref{fig:cube}) and add a seventh vertex, 
denoted $v$, that is connected to all the vertices of the octahedron.
The graph shown in Figure~\ref{fig:1222}
maps to the 1-skeleton of the octahedron; each vertex is shown as a small
circle labeled by the vertex of the octahedron it gets sent to.  It can be
checked directly that it is a finite planar emulator of the 1-skeleton
of the octahedron.

Note that some of the faces have been shaded (this includes the
outside face).  We add a vertex in each of these faces.  These
vertices all map to $v$ and are connected to every vertex on the
boundary of the shaded cells.  On the
boundary of each shaded face we see all the labels, so each of the
vertices that map to $v$ has all the necessary neighbors.  Finally,
we see that each vertex in Figure~\ref{fig:1222} is on the
boundary of at least one shaded face; hence, every vertex has a
neighbor that maps to $v$.   

This completes our construction of a finite planar emulator of
$K_{1,2,2,2}$. 

\medskip

\noindent{\it Remark.}  
By viewing $S^2$ as the boundary of the convex hull of the midpoints
of the edges of the cube or the octahedron, we may  draw the emulator
for $K_{1,2,2,2}$ symmetrically, so that it is invariant under the
antipodal involution.   The quotient gives an emulator for
$K_{1,2,2,2}$ that has 133 vertices and embeds in
$\mathbb{R}P^2$.  This symmetric presentation of the
emulator of $K_{1,2,2,2}$ reveals another interesting property.  By
considering the eight triangular faces (each shown in Figure~\ref{fig:1222} 
as a white triangle with a single shaded face), we can see
that they are formed from the union of 4 great circles, 
one with the labels 0, 1, 2, one with the labels 2, 3, 4,
one with the labels 1, 3, 5, and one with the labels 0, 4, 5.  Note
that if we two color the faces of the octahedron, we have exactly all
the faces of one color.


\begin{thebibliography}{10}

\bibitem{archdeacon-projective}
Dan Archdeacon.
\newblock A {K}uratowski theorem for the projective plane.
\newblock {\em J. Graph Theory}, 5(3):243--246, 1981.

\bibitem{archdeacon2}
Dan Archdeacon.
\newblock Two graphs without planar covers.
\newblock {\em J. Graph Theory}, 41(4):318--326, 2002.

\bibitem{fellowsphd}
Michael Fellows.
\newblock {\em Encoding Graphs in Graphs}.
\newblock PhD thesis, Univ. of California, San Diego,, 1985.

\bibitem{fellows}
Michael Fellows.
\newblock Planar emulators and planar covers.
\newblock 1988.

\bibitem{103}
Henry~H. Glover, John~P. Huneke, and Chin~San Wang.
\newblock 103 graphs that are irreducible for the projective plane.
\newblock {\em J. Combin. Theory Ser. B}, 27(3):332--370, 1979.

\bibitem{hlineny}
Petr Hlin{\v{e}}n{\'y}.
\newblock A note on possible extensions of {N}egami's conjecture.
\newblock {\em J. Graph Theory}, 32(3):234--240, 1999.

\bibitem{hlineny2}
Petr Hlin{\v{e}}n{\'y}.
\newblock 20 years of Negami's planar cover conjecture.
\newblock In {\em 20th Workshop on topological graph theory in Yokohama}, pages
  50--59, 2008.

\bibitem{kitakubo}
Shigeru Kitakubo.
\newblock Planar branched coverings of graphs.
\newblock {\em Yokohama Math. J.}, 38(2):113--120, 1991.

\bibitem{mohar}
Bojan Mohar.
\newblock Projective planarity in linear time.
\newblock {\em J. Algorithms}, 15(3):482--502, 1993.

\bibitem{negami}
Seiya Negami.
\newblock The spherical genus and virtually planar graphs.
\newblock {\em Discrete Math.}, 70(2):159--168, 1988.

\bibitem{negami-survey}
Seiya Negami.
\newblock Topological graph theory from {J}apan.
\newblock In {\em Proceedings of the {W}orkshop on {G}raph {T}heory and
  {R}elated {T}opics ({S}endai, 1999)}, volume~7, pages 99--112, 2001.

\bibitem{yr}
Yo'av Rieck and Yasushi Yamashita.
\newblock On Negami's planar cover conjecture.
\newblock available at http://arxiv.org/abs/math.CO/0612342, 2006.

\bibitem{rs2}
Neil Robertson and P.~D. Seymour.
\newblock Graph minors. {XIII}. {T}he disjoint paths problem.
\newblock {\em J. Combin. Theory Ser. B}, 63(1):65--110, 1995.

\bibitem{rs}
Neil Robertson and P.~D. Seymour.
\newblock Graph minors. {XX}. {W}agner's conjecture.
\newblock {\em J. Combin. Theory Ser. B}, 92(2):325--357, 2004.

\end{thebibliography}
\end{document}